\documentclass[12pt,reqno]{amsart}
\usepackage{amssymb, amsfonts, amsbsy, latexsym}
\usepackage[english]{babel}

\newtheorem{Thm}{Theorem}[section]

\newtheorem{corollary}[Thm]{Corollary}
\newtheorem{lemma}[Thm]{Lemma}
\newtheorem{proposition}[Thm]{Proposition}

\newtheorem{theorem}[Thm]{Theorem}

\def \H{\mathbb H}

\def \beq{\begin{eqnarray*}}
\def \eeq{\end{eqnarray*}}

\def\span{{\text{\rm span}}}
\newcommand{\ip}[2]{\left\langle#1,#2\right\rangle}
\newcommand{\abs}[1]{\left| #1 \right|}

\begin{document}

\title{A fundamental identity for Parseval frames}
\author[R. Balan, P.G. Casazza, D. Edidin, and G. Kutyniok]{Radu Balan, 
Peter G. Casazza, Dan Edidin, and Gitta Kutyniok}

\address{\textrm{(R.~Balan)}
Siemens Corporate Research,
755 College Road East,
Princeton, NJ 08540}
\email{radu.balan@siemens.com}

\address{\textrm{(P.~G.~Casazza)}
Department of Mathematics,
University of Missouri,
Columbia, MO 65211}
\email{pete@math.missouri.edu}

\address{\textrm{(D.~Edidin)}
Department of Mathematics,
University of Missouri,
Columbia, MO 65211}
\email{edidin@math.missouri.edu}

\address{\textrm{(G.~Kutyniok)}
Mathematical Institute,
Justus-Liebig-University Giessen,
35392 Giessen,
Germany}
\email{gitta.kutyniok@math.uni-giessen.de}

\begin{abstract}
In this paper we establish a surprising fundamental identity for 
Parseval frames in a Hilbert space. Several variations of
this result are given, including an extension to general frames.
Finally, we discuss the derived results.
\end{abstract}



\date{\today}

\maketitle

\section{Introduction}\label{Intro}

Frames are an essential tool for many emerging applications such as
data transmission.  Their main advantage is the fact that frames can
be designed to be redundant while still providing reconstruction
formulas.  This makes them robust against noise and losses while
allowing freedom in design (see, for example, \cite{CK01,GKK01}).  Due
to their numerical stability, {\it tight frames} and
{\it Parseval frames} are of increasing
interest in applications (See Section 2.1 for definitions.).
Particularly in image processing, tight frames have emerged as
essential tool (compare \cite{CRSS04}). In abstract frame theory,
systems constituting tight frames and, in particular, Parseval frames
have already been extensively explored \cite{BF03,CK01,CK05,EF02,
GKK01,VW05}, yet many questions are still open.

For many years engineers believed that, in applications
such as speech recognition, a signal can be reconstructed without information
about the phase.  In \cite{BCE} this
longstanding conjecture was verified by 
constructing new classes of Parseval frames
for which a signal vector can reconstructed without noisy phase or its
estimation.  While working on efficient algorithms for signal
reconstruction, the authors of \cite{BCE} discovered a surprising
identity for Parseval frames (see \cite{BCEK} for a detailed
discussion of the origins of the identity).

Our Parseval frame identity can be stated as follows (Theorem \ref{PFI}):
For any Parseval frame $\{f_i\}_{i\in I}$ in a Hilbert space $\H$, and
for every subset $J\subset I$ 
and every $f\in \H$ 
\begin{equation}\label{eq1}
\sum_{i\in J}|\langle f,f_i \rangle |^2 -\|\sum_{i\in J}
\langle f,f_i \rangle f_i \|^2 =
\sum_{i\in J^c}|\langle f,f_i \rangle |^2 -\|\sum_{i\in J^c}
\langle f,f_i \rangle f_i \|^2.
\end{equation}
The proof given here, based on operator theory,
admits an elegant extension to arbitrary frames (Theorem \ref{TT}). 
However, our main focus will be on Parseval
frames because of their importance in applications, 
particularly to signal processing.
Several interesting variants of our result are presented; 
for example, we
show that overlapping divisions can be also used.  Then the Parseval frame
identity is discussed in detail; in particular, we derive intriguing
equivalent conditions for both sides of the identity to be equal to
zero.

%This paper is organized as follows. In Section \ref{sec_napr} we
%introduce some notion and state some preliminary results. The Parseval
%Frame Identity (Theorem \ref{PFI}) and several related versions as
%well as a version for general frames are established in Section
%\ref{sec_afi}. In the last section (Section \ref{sec_dotpfi}) the
%derived results are discussed and fit into the theory of frames.

%******************************************************************************

\section{Notation and preliminary results}
\label{sec_napr}

%*************************************************************************************

\subsection{Frames and Bessel sequences}

Throughout this paper $\H$ will always denote a Hilbert space and $I$
an indexing set.  The finite linear span of a sequence of elements
$\{f_i\}_{i \in I}$ of $\H$ will be denoted by $\span(\{f_i\}_{i \in
I})$. The closure in $\H$ of this set will be denoted by
$\overline{\span}(\{f_i\}_{i \in I})$.

A system $\{f_i\}_{i \in I}$ in $\H$ is called a \emph{frame} for
$\H$, if there exist $0 < A \le B < \infty$ (\emph{lower and upper
frame bounds}) such that
\[A \, \|f\|_2^2
\le \sum_{i \in I}
    \abs{\ip{f}{f_i}}^2 \le B \, \|f\|_2^2 \quad \mbox{for all } f \in \H.\]
If $A,B$ can be chosen such that $A=B$, then $\{f_i\}_{i \in I}$ is an
\emph{$A$-tight frame}, and if we can take $A=B=1$, it is called a
\emph{Parseval frame}.  A \emph{Bessel sequence} $\{f_i\}_{i \in I}$
is only required to fulfill the upper frame bound estimate but not
necessarily the lower estimate. And a sequence $\{f_i\}_{i \in I}$ is
called a {\em frame sequence}, if it is a frame only for
$\overline{\span}(\{f_i\}_{i \in I})$.

The {\em frame operator} $Sf = \sum_{i \in I} \ip{f}{f_i}f_i$
associated with $\{f_i\}_{i \in I}$ is a bounded, invertible, and
positive mapping of $\H$ onto itself. This provides the frame
decomposition
\[f = S^{-1}Sf = \sum_{i \in I} \ip{f}{f_i} \tilde{f_i} = \sum_{i \in I} 
\langle f,\tilde{f_i}\rangle f_i,\] 
where $\tilde{f_i} = S^{-1}f_i$. The family $\{\tilde{f_i}\}_{i \in
I}$ is also a frame for $\H$, called the {\em canonical dual frame} of
$\{f_i\}_{i \in I}$.  If $\{f_i\}_{i\in I}$ is a Bessel sequence in
$\H$, for every $J\subset I$ we define the operator $S_J$ by
\[S_{J}f = \sum_{i\in J}\langle f,f_i \rangle f_i.\]

Finally, we state a known result (see, for example, \cite{C}), since it will
be employed several times.

\begin{proposition}\label{PPP}
Let $\{f_i\}_{i\in I}$ be a frame for $\H$
with frame operator $S$.  For every $f\in \H$, we have
\begin{enumerate}
\item  $\|\sum_{i\in I}\langle f,f_i \rangle f_i \|^2 \le \|S\| \, \sum_{i\in I}|\langle f,f_i\rangle |^2$.
\vspace*{0.2cm}
\item  $\sum_{i\in I}|\langle f,f_i\rangle |^2 \le \|S^{-1}\| \, \|\sum_{i\in I}\langle f,f_i \rangle f_i \|^2$.
\end{enumerate}
Moreover, both these inequalities are best possible.
\end{proposition}

For more details on frame theory 
 we refer to the survey article \cite{Cas00} and the book \cite{C}.
 
%*************************************************************************************

\subsection{Operator Theory}

We first state a basic result from Operator Theory, which is
very useful for the proof of the fundamental identity.

\begin{proposition}\label{P1}
If $S,T$ are operators on $\H$ satisfying $S+T=I$, then $S-T = S^2 -T^2$.
\end{proposition}

\begin{proof}
We compute
\[S-T = S-(I-S) =  2S-I = S^2 - (I-2S+S^2) = S^2 - (I-S)^2
= S^2 - T^2.\]
\end{proof}

\begin{proposition}
Let $S,T$ be operators on $\H$ so that
$S+T=I$.  Then $S,T$ are self-adjoint if and only if
$S^{*}T$ is self-adjoint.
\end{proposition}

\begin{proof}
Suppose that $S^{*}T$ is self-adjoint. Then
\[S-T = (S^{*}+T^{*})(S-T) = S^{*}S + T^{*}S - S^{*}T - T^{*}T
= S^{*}S-T^{*}T.\]
This shows that $S-T$ is self-adjoint.  Since $S+T$ is self-adjoint
by hypothesis, it follows that
\[S = \tfrac{1}{2}(S+T + (S-T)) \quad \mbox{and} \quad 
T = \tfrac12(S+T-(S-T))\]
are self-adjoint. 

The converse is obvious.
\end{proof}

%*************************************************************************************

\section{A fundamental Identity}
\label{sec_afi}

%*************************************************************************************

\subsection{General frames}

We first study the situation of general frames in $\H$. 

\begin{theorem}\label{TT}
Let $\{f_i \}_{i\in I}$ be a frame for $\H$ with canonical dual frame 
$\{\tilde{f}_i\}_{i\in I}$.  Then for all $J\subset I$ and all $f\in \H$ we have
$$
\sum_{i\in J}|\langle f,f_i \rangle |^2 - 
\sum_{i\in I}|\langle S_J f,\tilde{f}_i \rangle |^2 =
\sum_{i\in J^c}
|\langle f,f_i \rangle |^2 
 -\sum_{i\in I}|\langle S_{J^c}f,\tilde{f}_i \rangle |^2.
$$
\end{theorem}

\begin{proof}
Let $S$ denote the frame operator for $\{f_i\}_{i\in I}$. Since $S = S_J + S_{J^c}$, it
follows that $I = S^{-1}S_J + S^{-1}S_{J^c}$.
Applying Proposition \ref{P1} to the two operators $S^{-1}S_J$ and $S^{-1}S_{J^c}$
yields
\begin{equation}\label{eq11}
S^{-1}S_{J} - S^{-1}S_J S^{-1}S_J = S^{-1}S_{J^c} - S^{-1}S_{J^c}S^{-1}S_{J^c}.
\end{equation}
Thus for every $f,g\in \H$ we obtain 
\begin{equation}
\langle S^{-1}S_J f,g \rangle -
\langle S^{-1}S_J S^{-1}S_J f,g \rangle =
\langle S_{J}f,S^{-1}g \rangle -
\langle S^{-1}S_{J}f,S_{J}S^{-1}g \rangle. \label{E}
\end{equation}
Now we choose $g$ to be $g = Sf$. Then we can continue the
equality (\ref{E}) in the following way:
\[= \langle S_{J}f,f\rangle - \langle S^{-1}S_J f,S_{J} f
\rangle = \sum_{i\in J}|\langle f,f_i \rangle |^2
- \sum_{i\in I}|\langle S_{J}f,\tilde{f}_i \rangle |^2.\]
Setting equality (\ref{E}) equal to the corresponding
equality for $J^c$ and using \eqref{eq11}, we finally get
\[\sum_{i\in J}|\langle f,f_i \rangle |^2
- \sum_{i\in I}|\langle S_{J}f,\tilde{f}_i \rangle |^2
= \sum_{i\in J^c}|\langle f,f_i \rangle |^2
- \sum_{i\in I}|\langle S_{J^c}f,\tilde{f}_i \rangle |^2.\]
\end{proof}

%*************************************************************************************

\subsection{Parseval Frames}

In the situation of Parseval frames the fundamental identity is of a special
form, which moreover enlightens the surprising nature of it.

\begin{theorem}[Parseval Frame Identity]\label{PFI}
Let $\{f_i\}_{i\in I}$ be a Parseval frame for $\H$. For every subset $J\subset I$ and every
$f\in \H$, we have
\[\sum_{i\in J}|\langle f,f_i \rangle |^2 -\|\sum_{i\in J}
\langle f,f_i \rangle f_i \|^2 =
\sum_{i\in J^c}|\langle f,f_i \rangle |^2 -\|\sum_{i\in J^c}
\langle f,f_i \rangle f_i \|^2.\]
\end{theorem}

\begin{proof} 
We wish to apply Theorem \ref{TT}. Let $\{\tilde{f}_i\}_{i\in I}$ denote the dual frame
of $\{f_i\}_{i\in I}$. Since $\{f_i\}_{i\in I}$ is a Parseval frame,
its frame operator equals the identity operator and hence $\tilde{f}_i = f_i$ for
all $i\in I$. Employing Theorem \ref{TT} and the fact that $\{f_i\}_{i\in I}$ is a Parseval frame
yields
\begin{eqnarray*}
\sum_{i\in J}|\langle f,f_i \rangle |^2 - 
\|\sum_{i\in J}\langle f,f_i \rangle f_i \|^2 &=&
\sum_{i\in J}|\langle f,f_i \rangle |^2 -
\|S_{J}f\|^2 \\
&=&
\sum_{i\in J}|\langle f,f_i \rangle |^2 -
\sum_{i\in I}|\langle S_{J}f,f_i \rangle |^2 \\
&=&
\sum_{i\in J}|\langle f,f_i \rangle |^2 -
\sum_{i\in I}|\langle S_J f,\tilde{f}_i \rangle |^2\\ 
&=&
\sum_{i\in J^c}|\langle f,\tilde{f}_i \rangle |^2 -
\sum_{i\in I}
|\langle S_{J^c}f,\tilde{f}_i \rangle |^2 \\
&=& 
\sum_{i\in J^c}|\langle f,f_i \rangle |^2 
  - \|S_{J^c}f \|^2 \\
&=&
\sum_{i\in J^c}|\langle f,f_i \rangle |^2 -
\|\sum_{i\in J^c}\langle f,f_i \rangle f_i \|^2.
\end{eqnarray*}
\end{proof}

Note that the terms in the Parseval Frame Identity are
always positive (see Proposition \ref{PPP}).

A version of the Parseval Frame Identity for overlapping divisions is
derived in the following result.

\begin{proposition}
Let $\{f_i\}_{i\in I}$ be a Parseval frame for $\H$.
For every $J\subset I$, every $E\subset J^c$, and
every $f\in \H$, we have
\[\|\sum_{i\in J\cup E}\langle f,f_i \rangle f_i \|^2
-\| \sum_{i\in J^c \backslash E}\langle f,f_i \rangle f_i \|^2 =
\|\sum_{i\in J}\langle f,f_i \rangle f_i \|^2 - \|\sum_{i\in J^c}
\langle f,f_i \rangle f_i \|^2 + 2 \sum_{i\in E}|\langle
f,f_i \rangle |^2.\]
\end{proposition}

\begin{proof}
Applying Theorem \ref{PFI} twice yields
\begin{eqnarray*}
\|\sum_{i\in J\cup E}\langle f,f_i \rangle f_i \|^2
-\|\sum_{i\in J^c \backslash E}\langle f,f_i \rangle
\|^2 &=& \sum_{i\in J\cup E}|\langle f,f_i \rangle |^2 -
\sum_{i\in J^c \backslash E}|\langle f,f_i \rangle |^2 \\
&=& \sum_{i\in J}|\langle f,f_i \rangle |^2 - \sum_{i\in J^c}
|\langle f,f_i \rangle |^2 + 2 \sum_{i\in E}|\langle f,f_i
\rangle |^2 \\
&=& \|\sum_{i\in J}\langle f,f_i \rangle f_i \|^2 -
\|\sum_{i\in J^c}\langle f,f_i \rangle f_i \|^2
+ 2 \sum_{i\in E}|\langle f,f_i \rangle |^2.
\end{eqnarray*}
\end{proof}

Since each $\lambda$-tight frame can be turned into a Parseval frame by
a change of scale, we obtain the following corollary.

\begin{corollary}
\label{coro_tight}
Let $\{f_i\}_{i\in I}$ be a ${\lambda}$-tight frame
for $\H$.  Then for every $J\subset I$ and
every $f\in \H$ we have
\[{\lambda}\sum_{i\in J}|\langle f,f_i \rangle |^2 
-\|\sum_{i\in J}\langle f,f_i \rangle f_i \|^2 
={\lambda}\sum_{i\in J^c}|\langle f,f_i \rangle |^2
-\|\sum_{i\in J^c}\langle f,f_i \rangle f_i \|^2.\]
\end{corollary}

\begin{proof}
If $\{f_i\}_{i\in I}$ is a ${\lambda}$-tight frame for $\H$, then
$\{\frac{1}{\sqrt{\lambda}}f_i\}_{i\in I}$ is a Parseval frame
for $\H$. Applying Theorem \ref{PFI} proves the result.  
\end{proof}

Furthermore, the identity in Theorem \ref{PFI} remains true even for
Parseval frame sequences.

\begin{corollary}
Let $\{f_i \}_{i\in I}$ be a Parseval frame sequence
for $\H$.  Then for every $J\subset I$ and every $f\in \H$ we have
\[\sum_{i\in J}|\langle f,f_i \rangle |^2 -\|\sum_{i\in J}
\langle f,f_i \rangle f_i \|^2 =
\sum_{i\in J^c}|\langle f,f_i \rangle |^2 -\|\sum_{i\in J^c}
\langle f,f_i \rangle f_i \|^2.\]
\end{corollary}

\begin{proof}
Let $P$ denote the orthogonal projection of $\H$ onto
$\span(\{f_i \}_{i\in I})$.
By Theorem \ref{PFI}, we have
\[\sum_{i\in J}|\langle Pf,f_i \rangle |^2 -\|\sum_{i\in J}
\langle Pf,f_i \rangle f_i \|^2 =
\sum_{i\in J^c}|\langle Pf,f_i \rangle |^2 -\|\sum_{i\in J^c}
\langle Pf,f_i \rangle f_i \|^2.\]
Since $\langle Pf,f_i \rangle = \langle f,Pf_i \rangle
= \langle f,f_i \rangle$ for all $i \in I$, the result
follows.
\end{proof}

%*************************************************************************************

\section{Discussion of the Parseval Frame Identity}
\label{sec_dotpfi}

The identity given in Theorem \ref{PFI} is quite surprising
in that the quantities on the two sides of the identity are
not comparable to one another in general.  For example,
if $J$ is the empty set, then the left-hand-side of this
identity is zero because
$$
\sum_{i\in J}|\langle f,f_i
\rangle |^2 = 0 = \|\sum_{i\in J}\langle f,f_i \rangle f_i
\|^2.
$$
The right-hand-side of this identity is also zero, but now because
$$
\sum_{i\in J}|\langle f,f_i
\rangle |^2 = \|f\|^2 = \|\sum_{i\in J}\langle f,f_i \rangle f_i
\|^2.
$$
Similarly, if $|J|=1$, then both terms on
the left-hand-side of this identity may be arbitrarily close to
zero, while the two terms on the right-hand-side of the
identity are nearly equal to $\|f\|^2$, and they are canceling
precisely enough to produce the identity.  

If $\{f_i \}_{i\in I}$ is a Parseval frame for $\H$, then
for every $J\subset I$ and every $f\in \H$ we have
$$
\|f\|^2 = \sum_{i\in J}|\langle f,f_i \rangle |^2
+ \sum_{i\in J^c}|\langle f,f_i \rangle |^2.
$$
Hence, one of the two terms on the right-hand-side
of the above equality is greater than or equal to
$\frac{1}{2}\|f\|^2$.  It follows from Theorem
\ref{PFI} that for every $J\subset I$ and every
$f\in \H$,
$$
\sum_{i\in J}|\langle f,f_i \rangle |^2 +
\|\sum_{i\in J^c}\langle f,f_i \rangle f_i \|^2
= \sum_{i\in J^c}|\langle f,f_i \rangle |^2 +
\|\sum_{i\in J}\langle f,f_i \rangle f_i \|^2
\ge \frac{1}{2}\|f\|^2.
$$
We will now see that actually the right-hand-side 
of this inequality is in fact much larger.

\begin{proposition}
If $\{f_i\}_{i\in I}$ is a Parseval frame for
$\H$, then for every $J\subset I$ and every
$f\in \H$ we have
\[\sum_{i\in J}|\langle f,f_i \rangle |^2 +
\|\sum_{i\in J^c}\langle f,f_i \rangle f_i \|^2
\ge \tfrac{3}{4}\|f\|^2.\]
\end{proposition}

\begin{proof}
Since
\[\|f\|^2 = \| S_{J}f + S_{J^c}f \|^2 \le \|S_{J}f\|^2 + \|S_{J^c}f\|^2 
+ 2\|S_{J}f\|\|S_{J^c}f\| \le 2 ( \|S_{J}f\|^2 + \|S_{J^c}f\|^2 ),\]
we obtain
\[\langle (S_{J}^{2} + S_{J^c}^{2})f,f \rangle =
\|S_{J}f \|^2 + \|S_{J^c}\|^2 \ge \tfrac{1}{2}
\|f\|^2 = \langle \tfrac{1}{2}Id (f),f \rangle,\]
where $Id$ denotes the identity operator on $\H$.
Since $S_{J}+S_{J^c} = Id$, it follows that
$S_{J}+S_{J^c}^{2} +S_{J^c} + S_{J}^{2} \ge
\tfrac{3}{2}Id$.
Applying Proposition \ref{P1} to $S = S_{J}$ and $T=S_{J^c}$ yields 
$S_{J}+S_{J^c}^{2} = S_{J^c}+S_{J}^{2}$.  Thus
\[2(S_{J}+S_{J^c}^{2}) = S_{J}+S_{J^c}^{2}
+ S_{J^c} + S_{J}^2 \ge \tfrac{3}{2}Id.\]
Finally, for every $f\in \H$ we have
\[\sum_{i\in J}|\langle f,f_i \rangle |^2
+ \|\sum_{i\in J^c}\langle f,f_i \rangle f_i \|^2
= \langle S_{J}f,f \rangle + \langle S_{J^c}f,S_{J^c} f
\rangle = \langle (S_J + S_{J^c}^{2} )f,f \rangle 
\ge \tfrac{3}{4}\|f\|^2.\]
\end{proof}

Let $\{f_i\}_{i\in I}$ be a ${\lambda}$-tight frame for $\H$.  Reformulating 
Corollary \ref{coro_tight} yields that for every $J\subset I$ and
every $f\in \H$ we have
\[{\lambda}\sum_{i\in J}|\langle f,f_i \rangle |^2 
- {\lambda}\sum_{i\in J^c}|\langle f,f_i \rangle |^2
= \|\sum_{i\in J}\langle f,f_i \rangle f_i \|^2 
-\|\sum_{i\in J^c}\langle f,f_i \rangle f_i \|^2.\]
We intend to study when both sides of this equality equal zero, which is closely
related to questions concerning extending a frame to a tight frame. The proof
of this result uses the next lemma as a main ingredient.

\begin{lemma}\label{LL}
Let $\{f_i\}_{i\in I}$ and $\{g_i\}_{i\in K}$ be Bessel
sequences in $\H$ with frame operators $S$ and $T$, respectively.  
If $S=T$, then 
\[\span (\{f_i \}_{i\in I}) = \span (\{g_i\}_{i\in K}).\]
\end{lemma}

\begin{proof}
For any $f\in \H$, we have
\[\sum_{i\in I}|\langle f,f_i \rangle |^2 =
\langle Sf,f\rangle = \langle Tf,f \rangle =
\sum_{i\in K}|\langle f,g_i \rangle |^2.\]
It follows that $f\perp f_i$ for all $i\in I$ if
and only if $f\perp g_i$ for all $i\in K$.
\end{proof}

It is well known that given a frame $\{f_i\}_{i\in I}$ for a
Hilbert space $\H$, there exists a sequence (and in fact there
are many such sequences) $\{g_i\}_{i\in K}$ so that
$\{f_i\}_{i\in I} \cup \{g_i\}_{i\in K}$ is a tight frame.
We will now see that, if we choose two different families to
extend $\{f_i\}_{i\in I}$ to a tight frame, then these new
families have several important properties in common.

\begin{proposition}
Let $\{f_i\}_{i\in I}$ be a frame for $\H$.  Assume that $\{f_i \}_{i\in I} \cup 
\{g_i \}_{i\in K}$ and $\{f_i \}_{i\in I} \cup \{h_i \}_{i\in L}$ are both
${\lambda}$-tight frames.  Then the following condition hold.
\begin{enumerate}
\item For every $f\in \H$, $\sum_{i\in K}|\langle f,g_i \rangle |^2 = 
\sum_{i\in L}|\langle f,h_i \rangle |^2$.
\vspace*{0.2cm}
\item  For every $f\in \H$, $\sum_{i\in K}\langle f,g_i \rangle g_i = 
\sum_{i\in L}\langle f,f_i \rangle f_i$.
\vspace*{0.2cm}
\item  $\span (\{g_i\}_{i\in K}) = \span (\{h_i\}_{i\in L})$.
\end{enumerate}
\end{proposition}

\begin{proof}
For all $g\in \H$, we have
\[ \sum_{i\in I}|\langle f,f_i \rangle |^2 +
\sum_{i\in K}|\langle f,g_i \rangle |^2 =
{\lambda}\|f\|^2  =  \sum_{i\in I}|\langle f,f_i \rangle |^2 +
\sum_{i\in K}|\langle f,h_i \rangle |^2.\]
This yields (i).

Similarly,
\[\sum_{i\in I}\langle f,f_i \rangle f_i +
\sum_{i\in K}\langle f,g_i \rangle g_i =
{\lambda}f =
\sum_{i\in I}\langle f,f_i \rangle f_i +
\sum_{i\in K}\langle f,h_i \rangle h_i,\]
which proves (ii).

Condition (iii) follows immediately from (ii) and Lemma \ref{LL}.
\end{proof}

In the next result we will derive many equivalent conditions for both sides of the 
Parseval Frame Identity (Theorem \ref{PFI}) to equal zero. For this, we first need 
a technical result concerning the operators $S_J ,S_{J^c}$.  

\begin{proposition}\label{PP}
Let $\{f_i\}_{i\in I}$ be a Parseval frame for $\H$.  For any $J\subset I$, $S_J S_{J^c}$ 
is a positive self-adjoint operator on $\H$ which satisfies 
\[S_{J} - S_{J}^{2} = S_{J}S_{J^c} \ge 0.\] 
\end{proposition}

\begin{proof}
By symmetry and Proposition \ref{P1}, 
$S_J S_{J^c}$ is a positive self-adjoint operator on $\H$.  

Since $\{f_i\}_{i\in I}$ is a Parseval frame, for every $J\subset I$ and every
$f\in \H$, applying Proposition \ref{PPP} yields 
\[\langle S_{J}^{2}f,f \rangle =
\langle S_{J}f,S_{J}f \rangle =
\|\sum_{i\in J}\langle f,f_i \rangle f_i \|^2 \le 
\sum_{i\in I}|\langle f,f_i \rangle |^2 = 
\langle S_{J}f,f \rangle.\]
This proves $S_{J}-S_{J}^{2} \ge 0$.  Finally,
\[S_{J} = S_{J}(S_{J}+S_{J^c}) = S_{J}^{2} +
S_{J}S_{J^c}.\]
\end{proof}

Note that for any positive operator
$T$ on a Hilbert space $\H$ and any $f\in \H$, $Tf = 0$
implies $\langle Tf , f \rangle = 0$.  The converse of
this is also true.  If $\langle Tf,f\rangle = 0$, then
by a simple calculation
\[\langle Tf ,f \rangle = \langle T^{1/2}f,T^{1/2}f\rangle
= \|T^{1/2}f\|^2 = 0.\]
So $T^{1/2}f = 0$, and hence $Tf = 0$.
Noting that
one side of the Parseval Frame
Identity is zero if and only if the other side
is, we are led to the following result.

\begin{theorem}
Let $\{f_i \}_{i\in I}$ be a Parseval frame for $\H$. For each $J\subset I$ and $f\in \H$, 
the following conditions are equivalent.
\begin{enumerate}
\item  
$\sum_{i\in J}|\langle f,f_i \rangle |^2 =
\| \sum_{i\in J}\langle f,f_i \rangle f_i \|^2$.
\vspace*{0.2cm}
\item  
$\sum_{i\in J^c}|\langle f,f_i \rangle |^2 =
\| \sum_{i\in J^c}\langle f,f_i \rangle f_i \|^2$.
\vspace*{0.2cm}
\item  
$\sum_{i\in J}\langle f,f_i \rangle f_i \perp
\sum_{i\in J^c}\langle f,f_i \rangle f_i$.
\vspace*{0.2cm}
\item  $f \perp S_{J}S_{J^c}f$.
\vspace*{0.2cm}
\item  $S_{J}f = S_{J}^{2}f$.
\vspace*{0.2cm}
\item  $S_{J}S_{J^c}f = 0$.
\end{enumerate}
\end{theorem}

\begin{proof}

(i) $\Leftrightarrow$ (ii):  This is follows immediately from Theorem \ref{PFI}.

(iii) $\Leftrightarrow$ (iv):  This is proven by the following equality:
\[\langle \sum_{i\in J}\langle f,f_i \rangle f_i ,
\sum_{i\in J^c}\langle f,f_i \rangle f_i \rangle =
\langle S_{J}f,S_{J^c}f \rangle = \langle f,S_{J}S_{J^c}f \rangle.\]

(v) $\Leftrightarrow$ (vi):   This follows from
\[S_{J}^{2}f = S_{J}(I-S_{J^c})f = S_{J}f -
S_{J}S_{J^c}f.\]

(i) $\Leftrightarrow$ (v):  
We have
\[\sum_{i\in J}|\langle f,f_i \rangle |^2 - \|\sum_{i\in J}
\langle f,f_i \rangle f_i \|^2 = \langle S_J f,f\rangle
- \langle S_J f,S_J f \rangle = \langle (S_J -S_J^2 )f,f\rangle.\]
By Proposition \ref{PP}, $S_{J}-S_J^2 \ge 0$.  Therefore the
right-hand side of the above equality is zero if and only
if $(S_{J}-S_J^2)f = 0$ by our discussion preceding the
proposition. 

(i) $\Rightarrow$ (iv):  By (ii), $\langle S_Jf,f\rangle = \langle S_Jf,S_Jf \rangle$.
Hence $\langle (S_J-S_J^2 )f,f\rangle = \langle S_JS_{J^c}f,f\rangle
=0$, which implies (iv).

(iv) $\Rightarrow$ (vi):  By Proposition \ref{PP}, we have
that $S_J S_{J^c} \ge 0$.  Thus $\langle S_{J}S_{J^c}f,f\rangle 
= 0$ if and only if $S_{J}S_{J^c}f = 0$ by the discussion 
preceding this proposition.
\end{proof}

\section*{Acknowledgments}

An announcement for this paper appeared in \cite{BCEK}.

The authors wish to thank Alex Petukhov
for interesting discussions concerning this paper.  Petukhov
also provided us with an alternative matrix proof of the
Parseval Frame Identity (Theorem \ref{PFI}).  We also thank Chris
Lennard for useful discussions.  Lennard also provided us
with an alternative proof of the 
Parseval Frame Identity obtained by
expanding both sides as infinite series and comparing
the outcome. 

The second author was supported by NSF DMS 0405376,
the third author was supported by NSA MDA 904-03-1-0040,
and the fourth author was supported by DFG research
fellowship KU 1446/5.

%*************************************************************************************

\end{document}